\documentclass[12pt, a4paper]{article}
\usepackage{amsmath, amsthm, latexsym, enumerate}
\usepackage{amssymb}
\usepackage{graphicx} 
\usepackage{fullpage}
\usepackage{calc}
\usepackage[pdftex,colorlinks=true,urlcolor=blue,a4paper]{hyperref}
\usepackage{mathtools, times}
\usepackage[latin1]{inputenc}
\usepackage{amsmath}
\usepackage{multirow}
\usepackage{graphicx}
\usepackage[all]{xy}
\usepackage{amsfonts}
\usepackage{amsthm}
\usepackage{amssymb}
\usepackage{float}
\usepackage{geometry}
\usepackage{fullpage}
\theoremstyle{plain}

\newtheorem{theorem}{Theorem}[section]

\newtheorem{remark}{Remark}[section]

\numberwithin{equation}{section}
\date{}
\begin{document}
\begin{center} {A note on Andrews-MacMahon theorem}\end{center}
%\medskip
\begin{center}
 Darlison Nyirenda %\footnote{Corresponding Author: Darlison.Nyirenda@wits.ac.za} 
%\footnote{712040@students.wits.ac.za}
 \vspace{0.5cm} \\
 The John Knopfmacher Centre for Applicable Analysis and Number Theory, University of the Witwatersrand, P. O. Wits 2050, Johannesburg, South Africa.\\
e-mail: darlison.nyirenda@wits.ac.za\\
\end{center}
\begin{abstract}
\noindent For a positive integer $r$, George Andrews proved that the set of partitions of $n$  in which odd multiplicities are at least $2r + 1$ is equinumerous with the set of partitions of $n$ in which odd parts are congruent to $2r + 1$ modulo $4r + 2$. This was given as an extension of MacMahon's theorem ($r = 1$). Andrews, Ericksson, Petrov  and Romik gave a bijective proof of MacMahon's theorem. Despite several bijections being given, until recently, none of them was in the spirit of Andrews-Ericksson-Petrov-Romik bijection. Andrews' theorem has also been extended recently. Our goal is to give a generalized bijective mapping  of this further extension  in the spirit of Andrews-Ericksson-Petrov-Romik bijection. 
\end{abstract}

\section{Introduction}
A partition of $n$ is a representation $\lambda = \lambda_1 + \lambda_2 +  \cdots + \lambda_s$ in which the summands $\lambda_i$'s are integers such that $\sum\limits_{i = 1}^{s}\lambda_{i} = n$ and  $\lambda_1 \geq \lambda_{2} \geq \cdots \geq \lambda_{s}$. The summands are called \textit{parts} of $\lambda$ and the total sum $n$ is called the \textit{weight} of $\lambda$, denoted by $|\lambda|$.  A more convinient way of presenting $\lambda$ is to write
$\lambda = (\mu_1^{m_1}, \mu_{2}^{m_2}, \ldots \mu_{\ell}^{m_{\ell}})$
where $m_i$ is the multiplicity of the part $\mu_i$ (the number of times $\mu_i$ appears) and $\mu_1 > \mu_2 > \cdots > \mu_{\ell}$. Instead of considering any partition of fixed weight, sometimes, restrictions are imposed on the multiplicities as well as the parts themselves. When this happens, one may discover interesting identities. A classical example of this is the following theorem due to P. A. MacMahon.
\begin{theorem}[MacMahon, \cite{PAM}]\label{mac0}
The number of partitions of $n$ in which odd multiplicities are greater than 1 is equal to the number of partitions of $n$ in which odd parts are congruent to  $3 \pmod{6}$.
\end{theorem}
There are a couple of bijections for MacMahon's theorem and one of them  was given by Andrews, Ericksson, Petrov, and Romik \cite{gepr}. We shall call it the {\it Andrews-Ericksson-Petrov-Romik bijection} and describe it in the sequel.\\\\
\underline{Andrews-Ericksson-Petrov-Romik Bijection}\\\\
\noindent Let $A(n)$ denote the set of partitions of $n$ in which odd multiplicities are greater than 1 and let $B(n)$ denote the set of partitions of $n$ in which odd parts are congruent to $3 \pmod{6}$. \\
Let $ \lambda = (l^{h_{l}}, (l-1)^{h_{l-1}}, \ldots, 3^{h_{3}}, 2^{h_{2}}, 1^{h_{1}}) \in A(n)$ where $h_{i}$ is the multiplicity of $i$. Since no part appears exactly once, it follows that $h_{i} \in \{0, 2, 3, 4,\ldots \}$. Uniquely decompose $ h_{i}$ as $$h_{i} = k_{i} + g_{i},\,\,\text{where} \,\,\, k_{i} \in \{0, 3\},\,\, g_{i} \in \{0, 2, 4, 6, 8,\ldots \}.$$ 
For $j\geq 1$, define $d_j$ as  follows:
\begin{align*}
d_{6t + 1} & = d_{6t + 5} = 0, \\
d_{6t + 2} & = \frac{1}{2}g_{3t + 1}, \\
d_{6t + 4} & = \frac{1}{2}g_{3t + 2}, \\
d_{6t + 3} & = \frac{1}{3}k_{2t + 1} + g_{6t + 3}, \\
d_{6t + 6} & = \frac{1}{3}k_{2t + 2} + g_{6t + 6},
\end{align*} 
where $t = 0,1,2,\ldots$

The partition $(f^{d_f},(f-1)^{d_{f-1}},\ldots, 2^{d_2}, 1^{d_1})$ is in $B(n)$ and the transformation is invertible. 

\noindent Andrews gave the following generalisation of Theorem \ref{mac0}. 

\begin{theorem}[Andrews, \cite{andmac}] \label{mac}
Let $r \in \mathbb{Z}_{\ge 0}$. The set of partitions of $n$ in which parts with odd multiplicity appear at least $2r+1$ times is equinumerous with the set of partitions of $n$ in which parts are even or congruent to $2r+1 \pmod{4r+2}$.
\end{theorem}

Bijective proofs for Theorem \ref{mac} have been given (see \cite{JSSF} and the references therein). However, a bijection for this theorem that extends Andrews-Ericksson-Petrov-Romik  bijection has been given recently \cite{integer}.
In \cite{annals},  Theorem \ref{mac} was extended as follows.
\begin{theorem}[Nyirenda, Mugwangwavari, \cite{annals}]\label{nyirenda}
Let $a, p \in \mathbb{N}$ such that $\gcd(a,p) = 1$ and $a < p$. The number of partitions of $n$ in which multiplicities congruent to $ja \pmod{p}$ are at least $j(pr + a)$ for $j = 0,1,2\ldots, p-1$ is equal to the number of partitions of $n$ wherein parts not divisible by $p$ are congruent to $-s(pr + a) \pmod{p^2r + pa}$ where $s = 1,2, \ldots, p-1$. 
\end{theorem}
The authors, not only provided this extention of Andrews' theorem but also gave a bijection which naturally extends a bijection of Fu and Sellers \cite{annals}. However, finding a bijection for this extension that is in the spirit of the Andrews-Ericksson-Petrov-Romik  Bijection is still an open question. \\
In this paper we completely settle this question by providing a bijective proof of Theorem \ref{nyirenda} that is reminiscent of Andrews-Ericksson-Petrov-Romik bijection. The mapping we give also generalizes the bijection in \cite{integer}.
\section{The generalized bijection}
Let $\tilde{A}(n)$ be the set of partitions of $n$ in which multiplicities congruent to $ja \pmod{p}$ are at least $j(pr + a)$ for $j = 0,1,2\ldots, p-1$  and $\tilde{B}(n)$ denote the set of partitions of $n$  wherein parts not divisible by $p$ are congruent to $-s(pr + a) \pmod{p^2r + pa}$ where $s = 1,2, \ldots, p-1$. \\
\noindent Suppose that $ \lambda =  (l^{h_{l}}, (l-1)^{h_{l-1}}, \ldots, 3^{h_{3}}, 2^{h_{2}}, 1^{h_{1}}) \in \tilde{A}(n)$, Here, $h_{i}$ is the multiplicity of $i$. \\
We can decompose $ h_{i} $ uniquely as $h_{i} = k_{i} + g_{i}$ where $$ k_{i} = (pr + a)v_i \,\,\text{and}\,\, g_{i} = h_{i} - (pr + a)v_{i} .$$
The integer $v_i$ is the least positive residue of $a^{-1}h_i$ modulo $p$. Here $a^{-1}$ is the inverse of $a$ modulo $p$ (it exists because $\gcd(a,p) = 1$).
\noindent Since $0\leq v_i \leq  p - 1$, we must have $$ k_{i} \in \{0,pr+a,2(pr+a),\ldots,(p-1)(pr+a)\} $$ and $$ g_{i} \in \{0, p, 2p, 3p, 4p, 5p,\hdots \}.$$
\noindent  Define $d_j$ as  follows: 
\begin{align*}
d_{(p^2 r + pa)t + pj - a}  & =  0, j \in \{1,2,3,\ldots, pr + a \}\setminus \{ (p-1)r + a \} \\\\
d_{(p^2 r + pa)t + pj - 2a}  & =  0, j \in \{1,2,3,\ldots, pr + a \}\setminus \{ (p-2)r + a \} \\\\
d_{(p^2 r + pa)t + pj - 3a}  & =  0, j \in \{1,2,3,\ldots, pr + a \}\setminus \{ (p-3)r + a \}\\\\
                                    &\quad  \vdots   \\
                                    & \quad \vdots   \\\\
d_{(p^2 r + pa)t + pj - (p-1)a}  & = 0, j \in \{1,2,3,\ldots, pr + a \}\setminus \{ r + a \},
\end{align*}\\
$$ d_{(p^2 r + pa)t + pj} = \frac{1}{p}g_{(pr + a)t + j},  j = 1,2,3,\ldots, pr + a - 1$$\\
and\\
$$ d_{(p^2 r + pa)t  + (p - j)(pr + a)} =    \frac{1}{pr+ a}k_{pt + p - j} + g_{(p^2 r + pa)t + (pr + a)(p - j )}, j = 1, 2, 3, \ldots, p$$\\
with $t = 0, 1,2 \ldots$.\\\\

\noindent Let $\mu = (\ldots , i^{d_{i}},( i -1)^{d_{i - 1}}, (i - 2)^{d_{i - 2}}, \ldots 2^{d_2}, 1^{d_1})$. We claim that $\mu \in \tilde{B}(n)$. To see this, observe that,  if $f$ is a positive integer such that $f =(p^2 r + pa)t + pj - ia$ for some $t \geq0$, and $1\leq i\leq p - 1$  with  $j \in \{0,1,2,3,\ldots, pr + a \}\setminus \{ (p-i)r + a \}$, then
$f \equiv -ia \pmod{p}$ so that  $f \not\equiv 0 \pmod{p}$.  Also, observe that
\begin{align*}
  (p^2 r + pa)t + p((p - i)r + a) - ia  & =  (p^2 r + pa)t + p^{2}r  - ipr + pa - ia \\
    & =  (p^2 r + pa)(t + 1)   - i(pr + a) \\
     & \equiv - i(pr + a) \pmod{p^2 r + pa}.
\end{align*}
\noindent Since $j \neq (p - i)r + a$, we conclude that $f$ is not congruent to $-s(pr + a) \pmod{p^2r + pa}$ with $1\leq s \leq p - 1$. Thus, $d_{f} = 0$ implies that such integers cannot appear as parts, which is part of the description of partitions in $\tilde{B}(n)$.\\
\noindent Note that integers congruent to $0 \pmod{p}$ are allowed to appear, and these integers are precisely
those of the form $$(p^2 r + pa)t + pj\,\,\,\text{where}\,\,\,j = 0, 1, 2, \ldots, pr + a - 1\,\,\,\text{and} \,\,\, t \geq 0.$$
By definition of $d_{j}$'s, one can check that
$$
d_{(p^2 r + pa)t + pj} 
= \begin{cases}
\frac{1}{p}g_{pr + a)t + j}, & \text{if $j = 1,2,\ldots, pr + a - 1$;}\\\\
\frac{1}{pr+ a}k_{pt} + g_{(p^2 r + pa)t}, &\text{if $j = 0$,}
\end{cases}
$$ which is not necessarily 0, i.e. the multiplicity of a positive integer congruent to 0 modulo $p$ may be 0 or not. 
 Indeed, $\mu \in \tilde{B}(n)$. \\
\noindent Define a mapping $\phi: \tilde{A}(n) \rightarrow \tilde{B}(n)$ by $  \lambda \mapsto \mu$. In the next section, we prove that $\phi$ is a bijection from $\tilde{A}(n)$ onto $\tilde{B}(n)$.
\section{Proof of the bijection}
We first show that $\phi$  preserves weight, i.e. $|\mu| = n$.  Now 
\begin{align*}
|\mu|  & = | (\ldots , i^{d_{i}},( i -1)^{d_{i - 1}}, (i - 2)^{d_{i - 2}}, \ldots 2^{d_2}, 1^{d_1})| \\
          & = \sum_{i = 1}^{\infty}id_{i}\\
          & = \sum_{j = 1}^{(p-1)r + a - 1}\sum_{t = 0}^{\infty}\left((p^2 r + pa)t + pj - a \right)d_{(p^2 r + pa)t + pj - a} \\\\
          & \quad  \quad + \,\,  \sum_{j = (p - 1)r + a + 1}^{pr + a}\sum_{t = 0}^{\infty}\left((p^2 r + pa)t + pj - a \right)d_{(p^2 r + pa)t + pj - a}  \\\\
         &  \quad  \quad + \,\, \sum_{j = 1}^{(p-2)r + a - 1}\sum_{t = 0}^{\infty}\left((p^2 r + pa)t + pj - 2a \right)d_{(p^2 r + pa)t + pj - 2a} \\\\
          & \quad  \quad + \,\,  \sum_{j = (p - 2)r + a + 1}^{pr + a}\sum_{t = 0}^{\infty}\left((p^2 r + pa)t + pj - 2a \right)d_{(p^2 r + pa)t + pj - 2a}  \\
          & \quad  \quad + \,\, \sum_{j = 1}^{(p-3)r + a - 1}\sum_{t = 0}^{\infty}\left((p^2 r + pa)t + pj - 3a \right)d_{(p^2 r + pa)t + pj - 3a} \\\\
          & \quad  \quad + \,\,  \sum_{j = (p - 3)r + a + 1}^{pr + a}\sum_{t = 0}^{\infty}\left((p^2 r + pa)t + pj - 3a \right)d_{(p^2 r + pa)t + pj - 3a} \\\\
         & \quad  \quad + \,\, \sum_{j = 1}^{(p-4)r + a - 1}\sum_{t = 0}^{\infty}\left((p^2 r + pa)t + pj - 4a \right)d_{(p^2 r + pa)t + pj - 4a} \\
          & \quad  \quad + \,\,  \sum_{j = (p - 4)r + a + 1}^{pr + a}\sum_{t = 0}^{\infty}\left((p^2 r + pa)t + pj - 4a \right)d_{(p^2 r + pa)t + pj - 4a}  \end{align*}
\begin{align*}
          & \quad \quad \quad \qquad \vdots \\\\
           & \quad  \quad + \,\, \sum_{j = 1}^{r + a - 1}\sum_{t = 0}^{\infty}\left((p^2 r + pa)t + pj - (p - 1)a \right)d_{(p^2 r + pa)t + pj - (p-1)a} \\\\
          & \quad  \quad + \,\,  \sum_{j = r + a + 1}^{pr + a}\sum_{t = 0}^{\infty}\left((p^2 r + pa)t + pj - (p-1)a \right)d_{(p^2 r + pa)t + pj - (p-1)a} \\\\
          & \quad  \quad + \,\, \sum_{j = 1}^{pr + a - 1}\sum_{t = 0}^{\infty}\left((p^2 r + pa)t + pj\right)d_{(p^2 r + pa)t + pj} \\\\
          & \quad  \quad + \,\,   \sum_{j = 1}^{p}\sum_{t = 0}^{\infty} \left((p^2 r + pa)t  + (p - j)(pr + a) \right)d_{(p^2 r + pa)t  + (p - j)(pr + a)}\\
        & = \sum_{i = 1}^{p - 1} \sum_{j = 1}^{(p-i)r + a - 1}\sum_{t = 0}^{\infty}\left((p^2 r + pa)t + pj - ia \right)d_{(p^2 r + pa)t + pj - ia} \end{align*}
\begin{align*}
          & \quad  \quad + \,\,  \sum_{i = 1}^{p - 1}\sum_{j = (p - i)r + a + 1}^{pr + a}\sum_{t = 0}^{\infty}\left((p^2 r + pa)t + pj - ia \right)d_{(p^2 r + pa)t + pj - ia} \\\\
         & \quad  \quad + \,\,  \sum_{j = 1}^{pr + a - 1}\sum_{t = 0}^{\infty}\left((p^2 r + pa)t + pj\right)d_{(p^2 r + pa)t + pj} \\
          & \quad  \quad + \,\,  \sum_{j = 1}^{p}\sum_{t = 0}^{\infty} \left((p^2 r + pa)t  + (p - j)(pr + a) \right)d_{(p^2 r + pa)t  + (p - j)(pr + a)} \\\\
         & = \sum_{j = 1}^{pr + a - 1}\sum_{t = 0}^{\infty}\left((p^2 r + pa)t + pj\right)d_{(p^2 r + pa)t + pj} \\\\
          & \quad  \quad + \,\,   \sum_{j = 1}^{p}\sum_{t = 0}^{\infty} \left((p^2 r + pa)t  + (p - j)(pr + a) \right)d_{(p^2 r + pa)t  + (p - j)(pr + a)}\end{align*}
\begin{align*}
         &  = \sum_{j = 1}^{pr + a - 1}\sum_{t = 0}^{\infty}\left((p^2 r + pa)t + pj\right) \frac{1}{p}g_{pr + a)t + j} \\\\
          & \quad  \quad + \,\,  \sum_{j = 1}^{p}\sum_{t = 0}^{\infty} \left((p^2 r + pa)t  + (p - j)(pr + a) \right) \left(\frac{1}{pr+ a}k_{pt + p - j} + g_{(p^2 r + pa)t + (pr + a)(p - j )}\right)\\\\
         & = \sum_{j = 1}^{pr + a - 1}\sum_{t = 0}^{\infty}\left((pr + a)t + j\right)g_{(pr + a)t + j} +   \sum_{j = 1}^{p}\sum_{t = 0}^{\infty} \left(pt  + p - j \right)k_{pt + p - j} \\\\
        & \quad  \quad + \,\, \sum_{j = 1}^{p}\sum_{t = 0}^{\infty} \left((p^2 r + pa)t  + (p - j)(pr + a) \right) g_{(p^2 r + pa)t + (pr + a)(p - j )}.
\end{align*}
For simplification purposes, define summation operators $G$ and $H$ as follows:
$$ G(x,y) = \sum_{t = 0}^{\infty}\left(xt + y\right)g_{xt + y}$$
and 

$$ K(x,y) = \sum_{t = 0}^{\infty}\left(xt + y\right)k_{xt + y}.$$
\noindent Note that

\begin{align*}
\sum_{t = 0}^{\infty}\left((pr + a)t + j\right)g_{pr + a)t + j} & = \sum_{\ell = 0}^{p - 1}\sum_{t = 0}^{\infty}\left((p^2r + pa)t + j + (pr + a)\ell \right)g_{(p^2r + pa)t + j + (pr + a)\ell} \\
                                                                                          & = \sum_{\ell = 0}^{p - 1}G(p^2r + pa , j+ (pr + a)\ell),
\end{align*}
\begin{align*}
\sum_{t = 0}^{\infty} \left(pt  + p - j \right)k_{pt + p - j} & = \sum_{\ell = 0}^{pr + a - 1}\sum_{t = 0}^{\infty} \left((p^2r + pa)t  + p - j  + p\ell\right)k_{(p^2r + pa)t  + p - j  + p\ell} \\
                                                                                     & = \sum_{\ell = 0}^{pr + a - 1}K(p^2r + pa, p - j  + p\ell)
\end{align*}
and
$$ \sum_{t = 0}^{\infty} \left((p^2 r + pa)t  + (p - j)(pr + a) \right) g_{(p^2 r + pa)t + (pr + a)(p - j )} = G(p^2r + pa, (pr + a)(p - j )).$$
Hence,
\begin{align*}
|\mu|  & =  \sum_{j = 1}^{pr + a - 1} \sum_{\ell = 0}^{p - 1}G(p^2r + pa , j + (pr + a)\ell)   +   \sum_{j = 1}^{p}\sum_{\ell = 0}^{pr + a - 1}K(p^2r + pa, p - j  + p\ell)   \\\\
        & \quad  \quad + \,\, \sum_{j = 1}^{p}G(p^2r + pa, (pr + a)(p - j ))\\\\
        & =  \sum_{j = 1}^{pr + a - 1} \sum_{\ell = 0}^{p - 1}G(p^2r + pa , j + (pr + a)\ell)   +   \sum_{j = 0}^{p - 1}\sum_{\ell = 0}^{pr + a - 1}K(p^2r + pa, j  + p\ell)   \\\\
        & \quad  \quad + \,\, \sum_{\ell = 0}^{p - 1}G(p^2r + pa, (pr + a)\ell)\\\\
       &  =  \sum_{j = 0}^{pr + a - 1} \sum_{\ell = 0}^{p - 1}G(p^2r + pa , j + (pr + a)\ell)   +   \sum_{j = 0}^{p - 1}\sum_{\ell = 0}^{pr + a - 1}K(p^2r + pa, j  + p\ell) 
\end{align*}
\noindent Since modulo $pr + a$, the residues are $0,1,2, \ldots pr + a  - 1$, we have 
\begin{align*}
& \{j + (pr + a)\ell :   0\leq j \leq pr + a - 1,   0\leq \ell \leq p - 1\} \\
 & = \{ 0,1,2, \ldots, pr + a - 1 + (pr + a)(p-1) \} \\
                                                                                                & = \{0,1,2,3,\ldots, p^{2}r + pa - 1\}.\end{align*}

\noindent Similarly, modulo $p$, the residues are $0,1,2, \ldots, p-1$ and so we have
\begin{align*} \{j + p\ell :   0\leq j \leq p - 1,   0\leq \ell \leq pr + a - 1\}  & = \{ 0,1,2, \ldots, p - 1 + p(pr + a - 1) \}\\
 &  = \{0,1,2,3,\ldots, p^{2}r + pa - 1\}.\end{align*}
Thus,
$$ \{j + (pr + a)\ell :   0\leq j \leq pr + a - 1,   0\leq \ell \leq p - 1\} = \{j + p\ell :   0\leq j \leq p - 1,   0\leq \ell \leq pr + a - 1\}$$  which implies that
\begin{align*}
|\mu|  & = \sum_{j = 0}^{p^{2}r + pa - 1}G(p^2r + pa , j)   +   \sum_{j = 0}^{p^2 r + pa - 1}K(p^2r + pa, j)  \\
          & =  \sum_{j = 0}^{p^{2}r + pa - 1} \left(G(p^2r + pa , j)  + K(p^2r + pa, j)\right) \\
          & =   \sum_{j = 0}^{p^{2}r + pa - 1}\left( \sum_{t = 0}^{\infty}(p^2r + pa)t + j)g_{(p^2r + pa)t + j} + \sum_{t = 0}^{\infty}((p^2r + pa)t + j)k_{ (p^2r + pa)t  + j }\right) \\
        & = \sum_{j = 0}^{p^{2}r + pa - 1}\sum_{t = 0}^{\infty}\left( (p^2r + pa)t + j)g_{(p^2r + pa)t + j} + ((p^2r + pa)t + j)k_{ (p^2r + pa)t  + j }\right)\\
        & =  \sum_{j = 0}^{p^{2}r + pa - 1}\sum_{t = 0}^{\infty}(p^2r + pa)t + j)h_{(p^2r + pa)t + j)}\\
           & = \sum_{i = 0}^{\infty}ih_i\\
           & = |\lambda|.
\end{align*}

\noindent Finally, the transformation $\phi$ is reversible and its inverse is given as follows:\\
Let $\mu =  (\ldots , i^{d_{i}},( i -1)^{d_{i - 1}}, (i - 2)^{d_{i - 2}}, \ldots 2^{d_2}, 1^{d_1})$. For $j = 1, 2, \ldots,  pr + a -1$, we define
$$g_{(pr + a)t + j} = pd_{(p^2r + pa)t + j }, t \geq 0$$
and for $j = 0$, i.e. $g_{(pr + a)t}$, we  use:\\\\
For $j = 1,2,3,\ldots, p$,

$$
g_{(p^2r + pa)t + (pr + a)(p - j)} = 
\begin{cases}
d_{(p^2r + pa)t + (pr + a)(p - j)}, & \text{if $  d_{(p^2r + pa)t + (pr + a)(p - j)} \equiv 0 \pmod{p}$};\\\\
d_{(p^2r + pa)t + (pr + a)(p - j)} - 1, &  \text{if $  d_{(p^2r + pa)t + (pr + a)(p - j)} \equiv 1 \pmod{p}$};\\\\
d_{(p^2r + pa)t + (pr + a)(p - j)} - 2, & \text{if $  d_{(p^2r + pa)t + (pr + a)(p - j)} \equiv 2 \pmod{p}$};\\\\
d_{(p^2r + pa)t + (pr + a)(p - j)} - 3, &  \text{if $  d_{(p^2r + pa)t + (pr + a)(p - j)} \equiv 3 \pmod{p}$};\\\\
\quad \qquad \qquad \vdots  & \quad \qquad \qquad \vdots \\
d_{(p^2r + pa)t + (pr + a)(p - j)} - (p - 1), &  \text{if $  d_{(p^2r + pa)t + (pr + a)(p - j)} \equiv - 1 \pmod{p}$}
\end{cases}
$$

and
$$
k_{pt + p - j} = 
\begin{cases}
0, & \text{if $ d_{(p^2r + pa)t + (pr + a)(p - j)} \equiv 0 \pmod{p}$};\\\\
pr + a, &  \text{if $  d_{(p^2r + pa)t + (pr + a)(p - j)} \equiv 1 \pmod{p}$};\\\\
2(pr + a), & \text{if $ d_{(p^2r + pa)t + (pr + a)(p - j)} \equiv 2 \pmod{p}$};\\\\
3(pr + a), &  \text{if $  d_{(p^2r + pa)t + (pr + a)(p - j)} \equiv 3 \pmod{p}$};\\\\
\quad \qquad \vdots  & \quad \qquad \qquad \vdots \\
(p - 1)(pr + a), &  \text{if $  d_{(p^2r + pa)t + (pr + a)(p - j)} \equiv  p - 1 \pmod{p}$}.
\end{cases}
$$
Then the partition
$$\lambda = (\ldots, v^{k_{v} + g_v}, (v - 1)^{k_{v - 1} + g_{v - 1}}, (v - 2)^{k_{v - 2} + g_{v - 2} }, \ldots, 2^{k_2 + g_2}, 1^{k_1 + g_1})$$ is  in $\tilde{A}(n)$.
\begin{remark}
\noindent Setting $r = 1, a = 1, p = 2$ in the bijection yields the Andrews-Eriksson-Petrov-Romik bijction \cite{gepr}, $p = 2, a = 1$ yields Mugwangwavari-Nyirenda map  for general $r$ \cite{integer}.
\end{remark}
%\noindent For example, consider $n = 15, r = 3, a = 2$ and $p = 3$. Table \ref{tab1} demonstrates the correspondence.
%\begin{table}[h!]
%	\begin{center}
%		\begin{tabular}{ccc}
%			$A_{2}(15)$ & $\stackrel{}{\longrightarrow}$& $B_{2}(15)$\\ \hline
%			(2,2,2,2,2,1,1,1,1,1) & $\mapsto$ & (10,5) \\
%			(4,4,1,1,1,1,1,1,1) & $\mapsto$ & (8,5,2)\\
%			(3,3,2,Correct these,1,1) & $\mapsto$ & (6,5,4)\\
%			(3,3,1,1,1,1,1,1,1,1,1) & $\mapsto$ & (6,5,2,2)\\
%			(2,2,2,2,1,1,1,1,1,1,1) & $\mapsto$ & (5,4,4,2)\\
%			(2,2,1,1,1,1,1,1,1,1,1,1,1) & $\mapsto$ & (5,4,2,2,2)\\ 
%			(1,1,1,1,1,1,1,1,1,1,1,1,1,1,1) & $\mapsto$ & (5,2,2,2,2,2)	\\ \hline
%		\end{tabular}
%		\caption{The map $\phi: \tilde{A}(n)\rightarrow \tilde{B}(n)$ for $p = r = 3, a = 2, n=15$.}\label{tab1}
%	\end{center}
%\end{table}

\end{document}